\documentclass[aop,MSNbibl,number,dvips]{arximspdf}
\usepackage{graphicx}
\doi{10.1214/11-AOP718} 
\volume{41}
\issue{3B}
\pubyear{2013}
\firstpage{2261}
\lastpage{2278}

\makeatletter
\newtheorem{theorem}{Theorem}

\newtheorem{proposition}{Proposition}

\newproclaim{remark}{Remark}

\newcommand{\MaxFacets}{\operatorname{MaxPolytopes}_{d-1}}
\newcommand{\Cells}{\operatorname{Cells}}
\newcommand{\TypicalCell}{\operatorname{TypicalCell}}

\newcommand{\Vol}{\operatorname{Vol}}

\newcommand{\Var}{\operatorname{Var}}

\newcommand{\PHT}{\operatorname{PHT}}

\newcommand{\diam}{\operatorname{diam}}

\makeatother

\begin{document}
\begin{frontmatter}

\title{Limit theorems for iteration stable tessellations}
\runtitle{Limit theorems for STIT tessellations}

\begin{aug}
\author[A]{\fnms{Tomasz} \snm{Schreiber}\thanksref{t1}}
\and
\author[B]{\fnms{Christoph} \snm{Th\"ale}\corref{}\ead[label=e2]{cthaele@uos.de}}
\runauthor{T. Schreiber and C. Th\"ale}
\thankstext{t1}{Born June 25, 1975; died on December 1, 2010.}
\affiliation{Nicolaus Copernicus University Toru\'n and University of
Osnabr\"uck}
\address[A]{Faculty of Mathematics \\
\quad and Computer Science\\
Nicolaus Copernicus University\\
Toru\'n\\
Poland} 
\address[B]{Institut f\"ur Mathematik\\
Universit\"at Osnabr\"uck\\
Albrechtstra\ss e 28a\\
D-49076 Osnabr\"uck\\
Germany\\
\printead{e2}}
\end{aug}

\received{\smonth{4} \syear{2011}}
\revised{\smonth{9} \syear{2011}}

%
\begin{abstract}
The intent of this paper is to describe the large scale asymptotic
geometry of iteration stable (STIT) tessellations in $\mathbb{R}^d$,
which form a rather new, rich and flexible class of random
tessellations considered in stochastic geometry. For this purpose,
martingale tools are combined with second-order formulas proved earlier
to establish limit theorems for STIT tessellations. More precisely, a
Gaussian functional central limit theorem for the surface increment
process induced a by STIT tessellation relative to an initial time
moment is shown. As second main result, a central limit theorem for the
total edge length/facet surface is obtained, with a normal limit
distribution in the planar case and, most interestingly, with a
nonnormal limit showing up in all higher space dimensions.
\end{abstract}

%
\begin{keyword}[class=AMS]
\kwd[Primary ]{60D05}
\kwd{60F17}
\kwd[; secondary ]{60F05}
\kwd{60J75}.
\end{keyword}

\begin{keyword}
\kwd{Central limit theorem}
\kwd{functional limit theorem}
\kwd{iteration/nesting}
\kwd{Markov process}
\kwd{martingale theory}
\kwd{random tessellation}
\kwd{stochastic stability}
\kwd{stochastic geometry}.
\end{keyword}

\end{frontmatter}
%

\section{Introduction and results}\label{sec1}

Random tessellations or mosaics of $\mathbb{R}^d$ (with $d\geq2$) are
locally finite families of compact convex random polytopes, which have
no common interior points and cover the whole space. They form a
central object studied in stochastic geometry, spatial statistics and
related fields; see~\cite{SW,SKM} and the references cited therein.
However, there are only very few mathematically tractable models. The
most prominent examples include hyperplane and Voronoi tessellations,
where most often the Poisson case is considered. A new class, the
so-called STIT tessellations, was introduced recently by Mecke, Nagel
and Wei\ss\ in~\cite{MNW,MNW2,NW03,NW05} and has quickly attracted
considerable interest. These tessellations clearly show the potential
to become a new reference model for both theoretical and practical
purposes. Whereas most research on random tessellations in the last
decades has been about mean values and mean value relations (see \cite
{SW} for the recent state of the art), modern stochastic geometry
focuses on distributional aspects (\cite{BL2,HS1}, e.g.) and limit
theorems; see~\cite{HSS,HM} and the references therein. This paper
adds to these recent findings by providing a limit theory for STIT
tessellations.

In contrast to the tessellations studied so far, the STIT model has the
additional feature of arising as a result of a spatiotemporal \textit
{dynamic} construction. From this point of view, limit theorems for
STIT tessellations are particularly interesting. As recently pointed
out in~\cite{Schr}, we expect that the large scale asymptotic of
dynamic models for spatial random structures will become of great
importance in stochastic geometry and its applications in the near future.

\begin{figure}

\includegraphics{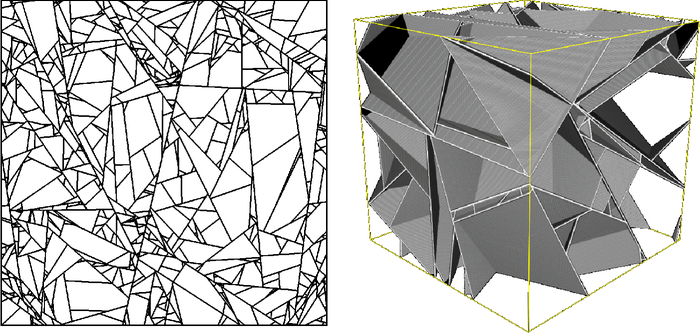}

\caption{Simulations of a planar and a spatial stationary and
isotropic STIT tessellation (kindly provided by Joachim Ohser and
Claudia Redenbach).}\label{Fig1}
\end{figure}

Let us recall the basic construction of tessellations that arise as a
result of repeated cell division. To this end, we identify the space
$\mathcal{ H}$ of hyperplanes in $\mathbb{R}^d$ with the parameter space $\mathbb{R}_+\times\mathcal{S}_{d-1}$
and the hyperplane $\{x\in\mathbb{R}^d:\langle
x,u\rangle=r\}$ with the pair $(r,u)\in\mathbb{R}_+\times\mathcal{S}_{d-1}$
and let $\Lambda$ be a measure on $\mathcal{H}$ which admits under the
described polar identification a representation of the form
%
\begin{equation}\label{LADEF}
\Lambda=\ell_+ \otimes\mathcal{R},
\end{equation}
where $\ell_+$ is the Lebesgue measure on the positive real half-axis
and where $\mathcal{R}$ is a probability measure on the unit sphere $\mathcal{S}_{d-1}$.
Throughout this paper we always require that the support of
$\mathcal{R}$ spans the whole space, that is, that $\operatorname{span}(\operatorname
{supp}(\mathcal{R}))=\mathbb{R}^d$, and we say in this case that $\Lambda$
is nondegenerate. Further, let $t>0$ be fixed, and let $W\subset\mathbb{R}^d$ be a
compact convex set with interior points in which our
construction of a random tessellation $Y(t\Lambda,W)$ is carried out.
In a first step, we assign to the window $W$ a random lifetime. Upon
expiry of its lifetime, the primordial cell $W$ dies and splits into
two sub-cells $W^+$ and $W^-$ separated by a hyperplane hitting~$W$,
which is chosen according to suitable restriction and normalization of~$\Lambda$.
The resulting new cells $W^+$ and $W^-$ are again assigned
independently of each other with random lifetimes and the entire
construction continues recursively until the previously fixed
deterministic time threshold $t$ is reached. The described process of
recursive cell divisions is called the \textit{MNW-construction}
through this paper (M-N-W stand for the inventors of the model), and
the resulting random subdivision of $W$ is denoted by $Y(t\Lambda,W)$;
see Figure~\ref{Fig1} for illustrations of the outcome of the
MNW-construction for $d=2$ and $d=3$. Note that the cells of $Y(t\Lambda
,W)$ are polyhedral except possibly those hitting the potentially
curved boundary of $W$, so that upon boundary effects $Y(t\Lambda,W)$
is a random tessellation of $W$.

In order to ensure the Markov property of the above construction in the
continuous-time parameter $t$, we assume from now on that the lifetimes
arising in the MNW-construction (including that of the initial window
$W$) are exponentially distributed. Moreover, we assume that the
parameter of the exponentially distributed lifetime of an individual
cell $[c]$ equals $\Lambda([c])$, where $[c]$ stands for the set of all
parameter values of hyperplanes hitting $c$. In this special situation,
the random tessellations $Y(t\Lambda,W)$ fulfill a stochastic stability
property under the operation of iteration of tessellations and are for
this reason called random \textit{STIT tessellations}; see Section \ref
{secSTIT} below for details.

%

Having studied the first- and second-order properties of STIT
tessellations in~\cite{STP1,STP2}, we consider in this paper the
central limit problem. This problem will be approached in two closely
related settings, interestingly leading to results of very different
qualitative nature. First, we shall focus our interest on the residual
length/surface increment process, arising, respectively, as cumulative
length or surface area of the cell-separating $(d-1)$-polyhedral facets
born \textit{after} a certain fixed initial time in the course of the
MNW-construction. In this set-up we shall establish a central limit
theorem with a Gaussian limiting variable. Next, we shall pass to the
\textit{total} length/surface, taking into account also the
segments/facets born at the very initial \textit{big bang} stages of the
MNW-construction, as descriptively termed in~\cite{MNW}. It turns out
that, whereas in dimension $2$ the Gaussian convergence is preserved,
this is no more the case for dimensions $3$ and higher, where
non-Gaussian limits arise. This apparently surprising phenomenon is in
fact due to the influence of the big bang phase in the MNW-construction
itself, which is negligible in the planar case, but turns out to be
crucial in higher dimensions.

We are now going to describe some of our limit theorems in more detail.
For a compact convex set $W$ as above, we put $W_R := R W$ for $R > 0$
and let $\Lambda$ be some fixed hyperplane measure as in (\ref{LADEF}).
Therefore and in order to simplify the notation we will write from now
on $Y(t,W)$ instead of $Y(t\Lambda,W)$. Our first limit theorem deals
with the total surface area\vadjust{\goodbreak} $\Vol_{d-1}(Y(t,W))$ of cell boundaries
induced by the MNW-construction within the time period $[s_0,1]$, where
$0<s_0<1$ is some positive initial time moment.
\begin{theorem}\label{limit1}
For each $0<s_0<1$, the random variable
\begin{eqnarray*}
&&\frac{1}{R^{d/2}}\bigl[\bigl({\Vol}_{d-1}(Y(1,W_R))-\mathbb{E}{\Vol
}_{d-1}(Y(1,W_R))\bigr)
\\
&&\qquad{} - \bigl({\Vol}_{d-1}(Y(s_0,W_R))-\mathbb{E}{\Vol}_{d-1}(Y(s_0,W_R))\bigr)\bigr]
\end{eqnarray*}
converges, as $R\rightarrow\infty$, in law to $\mathcal{N}(0,V_W(\Vol
_{d-1},\Lambda)\int_{s_0}^1 s^{1-d} \,ds)$, a normal distribution with
mean $0$ and variance $V_W(\Vol_{d-1},\Lambda)\int_{s_0}^1 s^{1-d} \,ds$,
where\break $V_W(\Vol_{d-1},\Lambda)$ is explicitly given by (\ref{VPHI}) or
alternatively (\ref{VPHIALTERNATIVE}) below.
\end{theorem}

This statement cannot be extended to $s_0 \downarrow0$, as would be of
interest as potentially leading to a Gaussian limit for the (centered
and suitably normalized) total edge length/surface area ${\Vol
}_{d-1}(Y(1,W_R)).$ The problem is that the variance integral $V_W(\Vol
_{d-1},\Lambda) \int_{s_0}^t s^{1-d} \,ds$ diverges at $0.$ However, this
difficulty can be overcome for $d=2$, but not for $d>2$. Indeed, in the
planar case the asymptotic behavior of the total edge length turns out
to be Gaussian:
\begin{theorem}\label{limit2}
For a STIT tessellation $Y(1,W_R)$ in the plane we have
\[
\frac{1}{R \sqrt{\log R}} [{\Vol}_1(Y(1,W_R))-\mathbb{E}{\Vol
}_1(Y(1,W_R))]\Longrightarrow\mathcal{N}(0,V_W(\Vol_1,\Lambda)),
\]
where $\Longrightarrow$ means convergence in law and where again
$V_W(\Vol_1,\Lambda)$ is explicitly known and given by (\ref{VPHI}) or
(\ref{VPHIALTERNATIVE}) below.
\end{theorem}
In fact, Theorems~\ref{limit1} and~\ref{limit2} are direct
consequences of our much stronger functional central limit theorems,
Theorems~\ref{BrownianConvergence} and~\ref{CLT2D} below.

For space dimensions $d>2$ we claim that the Gaussian convergence
cannot be preserved. Even though we are able to show this fact for all
$W$ and translation invariant $\Lambda$ by establishing non-Gaussian
tail decay, for simplicity, and in order to keep the argument
transparent, we only give a proof for a more easily tractable
particular case, in which all cells have the shape of cuboids
(rectangular parallelepipeds). The study of more involved properties of
the resulting random field is postponed to a future paper.
\begin{theorem}\label{thmNONGAUSSIAN} Fix $d>2$, take $W = [0,1]^d$ and
consider the hyperplane measure
%
\begin{equation}\label{PARTICULARLA}
\Lambda:= (2d)^{-1}\sum_{i=1}^d \int_{-\infty}^{+\infty} \delta_{r
e_i + e_i^{\bot}} \,dr,
\end{equation}
where $e_i,  i =1,\ldots,d$ are vectors of the standard orthonormal
basis for $\mathbb{R}^d$ and $\delta_{r e_i + e_i^{\bot}}$ is the unit
mass concentrated\vadjust{\goodbreak} on the hyperplane orthogonal to $e_i$ at distance $r$
from the origin (here $e_i^\bot$ stands for the orthogonal complement
of $e_i$). In this setting,
\[
{1\over R^{d-1}}[{\Vol}_{d-1}(1,W_R)-\mathbb{E}{\Vol}_{d-1}(1,W_R)]
\]
converges, as $R\rightarrow\infty$, to a non-Gaussian square-integrable
random variable $\Xi(W)$ with explicitly known variance given by (\ref
{XIVAR}) below.
\end{theorem}

The paper is structured as follows: In the next section we recall some
properties of STIT tessellations needed for the proofs of our limit
theorems. We also recall there some of the facts from~\cite{STP1,STP2}
in order to keep the paper self-contained and present the exact
statements of our functional central limit theorems. The proofs of our
results are the content of Section~\ref{secPROOFS}.

We would like to remark that an extended version~\cite{ST} of this text
is available online and, moreover, that the results in the present
paper form the basis of our works~\cite{ST2,ST3}.

\section{Background material and statement of the functional limit
theorems}\label{secSTIT}

We start by rephrasing some of the properties of the STIT tessellations
$Y(t\Lambda,W)$ as defined in the introduction, the proofs of which may
be found in~\cite{NW05}:
\begin{itemize}
\item$Y(t\Lambda,W)$ is consistent in that $Y(t\Lambda,W) \cap V
\stackrel{D}{=} Y(t\Lambda,V)$ for convex $V \subset W$ (here and below
$\stackrel{D}{=}$ stands for equality in distribution) and thus
$Y(t\Lambda,W)$ can be extended to a random tessellation $Y(t\Lambda)$
in the whole space $\mathbb{R}^d$ such that $Y(t\Lambda)\cap W$ has the
same distribution as $Y(t\Lambda,W)$.

This way, instead of interpreting $Y(t\Lambda,W)$ as outcome of the
MNW-construction carried out in $W$, $Y(t\Lambda,W)$ can also be
understood as restriction of the whole space random tessellation
$Y(t\Lambda)$ (which is a proper tessellation in the usual sense of
stochastic geometry as discussed in the introduction) to the window $W$.
\item$Y(t\Lambda)$ is a stationary random tessellation, that is,
stochastically translation invariant. If, moreover, $\Lambda$ is the
isometry-invariant hyperplane measure $\Lambda_\mathrm{iso}$, or
equivalently if $\mathcal{R}$ in (\ref{LADEF}) is the uniform distribution
$\nu_{d-1}$ on $\mathcal{S}_{d-1}$, then $Y(t\Lambda_\mathrm{iso})$ is even
isotropic, that is, stochastically invariant under rotations around the origin.
\item$Y(t\Lambda)$ is {st}able under the operation of {it}eration, denoted by $\boxplus$. This is to say
\[
Y(t\Lambda)\stackrel{D}{=}m\bigl(Y\bigl((t/m)\Lambda\bigr)\boxplus\cdots\boxplus
Y\bigl((t/m)\Lambda\bigr)\bigr)  \qquad   m=2,3,\ldots.
\]
For this reason, $Y(t\Lambda)$ is called a random STIT tessellation.
This property was discussed in detail in~\cite{NW05,STP1}, and we
refer to these papers and the references cited therein for further
discussion, because our arguments do not explicitly use the stochastic
stability, but its consequences.
\item The surface density of $Y(t\Lambda)$, that is, the mean surface
area of cell boundaries of $Y(t\Lambda)$ per unit volume equals $t$. In
particular, the mean surface area of facets arising in the
MNW-construction during time $[0,t]$ within a compact convex $W\subset
\mathbb{R}^d$ with interior points is given by $t\Vol_d(W)$.
\item STIT tessellations enjoy the following scaling property:
\[
tY(t\Lambda)\stackrel{D}{=}Y(\Lambda),
\]
that is, the
tessellation $Y(t\Lambda)$ of surface density $t$ upon rescaling by
factor $t$ has the same distribution as $Y(\Lambda)$, the STIT
tessellation with surface density $1$.
\item The intersection of a STIT tessellation in $\mathbb{R}^d$ with a
$k$-dimensional affine subspace $L_k$ ($1\leq k\leq d-1$), that is
independent of the tessellation induces a STIT tessellation in $L_k$.
\item STIT tessellations have Poisson typical cells, which is to say
that the distribution of the interior of the typical cell $\TypicalCell
(Y(t\Lambda))$ of $Y(t\Lambda)$ coincides with that of a Poisson
hyperplane tessellation $\PHT(t\Lambda)$ having intensity measure
$t\Lambda$; see the discussion after Theorem 1 in~\cite{STP1} or~\cite{NW03}.
\end{itemize}

For the nonspecialized reader let us remark that the typical cell of a
tessellation is what we get when we choose equiprobably a cell of the
tessellation at random out of a ``large'' observation window. The exact
definition makes use of Palm theory for which we refer to~\cite{SW}.
Moreover, let us recall that a Poisson hyperplane tessellation $\PHT
(t\Lambda)$ is a random subdivision of $\mathbb{R}^d$ induced by a
Poisson point process on the space of hyperplanes $\mathcal{H}$ with
intensity measure $t\Lambda$.

The finite volume continuous-time incremental MNW-construction of
random STIT tessellations, as discussed in the \hyperref[sec1]{Introduction}, clearly
has the Markov property in the continuous-time parameter, whence
natural martingales arise, which will be of crucial importance for our
further considerations. In fact, this observation was the starting
point of~\cite{STP1}, where a class of martingales associated to STIT
tessellations was constructed. In order to streamline our discussion we
do not repeat the full theory here, but rephrase the martingale
property of two stochastic processes, on which the proofs of our limit
theorems are based. To this end let $Y$ be some instant of $Y(t\Lambda
,W)$, and let $\phi(\cdot)$ be a measurable facet functional of the form
%
\begin{equation}
\phi(f) = \Vol_{d-1}(f) \zeta(\vec\mathbf{n}(f))\label{PHIFORM}
\end{equation}
with $\vec\mathbf{n}(f)$ standing for the unit normal to a facet $f$ and
$\zeta$ for a bounded measurable function on $\mathcal{S}_{d-1}$.
Moreover, denote the collection of cell-separating $(d-1)$-dimensional
facets, usually referred to as \textit{$(d-1)$-dimensional maximal
polytopes}, arising in subsequent splits in the MNW-construction by
$\MaxFacets(Y)$. (Note that some of these polytopes can be chopped-off
by the possibly curved boundary of the convex window $W$ in which $Y$
is constructed and are\vadjust{\goodbreak} no polytopes in the usual sense. However, we
somehow abuse notation and remark that this technical issue causes no
difficulties in our theory, because of the special form of the facet
functional $\phi$.) Let us further define $\Sigma_{\phi}(Y)$ by
\begin{equation}
\label{FDEF}
\Sigma_{\phi}(Y) := \sum_{f \in\MaxFacets(Y)} \phi(f)
\end{equation}
and $A_{\phi^2}(Y)$ to be
\[
A_{\phi^2}(Y) := \int_{[W]}\sum_{f\in\Cells(Y\cap H)}\phi^2(f)\Lambda(dH)
\]
with $\Cells(Y\cap H)$ standing for the set of $(d-1)$-dimensional
cells of the tessellation $Y\cap H$ induced by the intersection of $Y$
with a hyperplane $H$ (again, some of these cells may have a curved
boundary, because of the intersection with the construction window).
Let us further introduce the bar notation $\bar\Sigma_\phi(Y)$ for the
centred version $\Sigma_\phi(Y)-\mathbb{E}\Sigma_\phi(Y)$ of $\Sigma_\phi
(Y)$. Then, we have (see~\cite{STP1,STP2}):
\begin{proposition}\label{propMART3}
The two stochastic processes
%
\begin{equation}\label{MART3}
\bar\Sigma_{\phi}(Y(t\Lambda,W))\quad \mbox{and}\quad \bar\Sigma_{\phi
}^2(Y(t\Lambda,W)) - \int_0^t A_{\phi^2}(Y(s\Lambda,W))\,ds
\end{equation}
are both $\Im_t$-martingales, where $\Im_t$ stands for the filtration
generated by $(Y(s\Lambda,  W))_{0\leq s\leq t}$.
\end{proposition}

In particular (see \cite[Paragraph I.8]{LS}), by (\ref{MART3}), the
martingale $\bar\Sigma_{\phi}(Y(t\Lambda,W))$ has its
predictable quadratic variation process $\langle\bar\Sigma_{\phi
}(Y(\cdot,W)) \rangle$ absolutely continuous (in the sense of
functions) and
given by
%
\begin{equation}\label{PREDQUVAR}
\langle\bar\Sigma_{\phi}(Y(\cdot,W)) \rangle_t = \int_0^t A_{\phi
^2}(Y(s\Lambda,W)) \,ds.
\end{equation}

Besides these martingale tools, we will also make use of the following
formula for the variance $\Var(\Sigma_{\phi}(Y(t\Lambda,W)))$ of $\Sigma
_\phi(Y(t\Lambda,W))$, $W\subset\mathbb{R}^d$ compact, convex and with
interior points, established in full generality in~\cite{STP2}, in
order to calculate the variance of the limit random variable of our
non-Gaussian limit theorem:
\begin{proposition}\label{thmvariance} For any nondegenerate
translation-invariant $\Lambda$ of the form~(\ref{LADEF}) and $\phi$ as
in (\ref{PHIFORM}), we have
\begin{eqnarray*}
&&\Var\bigl(\Sigma_{\phi}(Y(t\Lambda,W))\bigr)\\
&&\qquad=\int_{[W]} \zeta
^2(\vec\mathbf{n}(H)) \int_{W \cap H} \int_{W \cap H}\frac{1-\exp
(-t\Lambda([xy]))}{\Lambda([xy])} \,dx \,dy \Lambda(dH),
\end{eqnarray*}
where $[xy]$ stands for the set of parameter values of hyperplanes
hitting the line segment $xy$ connecting $x$ with $y$.\vadjust{\goodbreak}
\end{proposition}
Let us further recall from~\cite{STP2} that the variance of the total
edge length of a stationary and isotropic STIT tessellation $Y(t\Lambda
_\mathrm{iso},W_R)$ in the plane behaves asymptotically like
%
\begin{equation}\label{EQVARASYM2D}
\Var\bigl(\Vol_{1}(Y(t\Lambda_\mathrm{iso},W_R))\bigr) \sim\pi\Vol_2(W)R^2\log R  \qquad  \mbox{as } R\rightarrow\infty,
\end{equation}
where $W$ is again a compact convex set as above. Indeed, this can be
seen from the general statement in Proposition~\ref{thmvariance}
combined with tools from integral geometry. Note that the asymptotic
variance expression for the total edge length is independent of $t$.
However, for all space dimensions $>2$, the surface density $t$ enters
the asymptotic variance expression as shown in detail in~\cite{STP2}.

We can now turn to the statement of our main results, the functional
limit theorems, from which Theorems~\ref{limit1} and~\ref{limit2} are
direct consequences. As in the \hyperref[sec1]{Introduction} we fix a hyperplane measure
$\Lambda$ and suppress from now on the reference to $\Lambda$, for
example by writing $Y(t,W_R)$ instead of $Y(t\Lambda,W_R)$.
\begin{theorem}\label{BrownianConvergence}
For each $0<s_0<1$ the centred surface increment process
\[
\biggl( \mathcal{S}^{R,W}_{s_0,t} :=
\frac{1}{R^{d/2}} [\bar\Sigma_{\phi}(Y(t,W_R)) - \bar\Sigma_{\phi
}(Y(s_0,W_R))] \biggr)_{t\in[s_0,1]}
\]
converges in law, as $R\to\infty,$ on the space $\mathcal{D}[s_0,1]$ of
right continuous functions with left-hand limits
(c\`adl\`ag) on $[s_0,1]$, endowed with
the Skorokhod $J_1$-topology \cite[Chapter 6.1]{LS}, to a time-changed
Wiener process
\[
t \mapsto\mathcal{W}_{V_W(\phi,\Lambda) \int_{s_0}^t s^{1-d} \,ds},
\]
where $\mathcal{W}_{(\cdot)}$ is the standard Wiener process, and $V_W(\phi
,\Lambda)$ is given
by (\ref{VPHI}), or alternatively (\ref{VPHIALTERNATIVE}), below. In particular,
\[
\mathcal{S}^{R,W}_{s_0,1} = \frac{1}{R^{d/2}} [\bar\Sigma_{\phi
}(Y(1,W_R)) - \bar\Sigma_{\phi}(Y(s_0,W_R))]
\]
converges in law to $\mathcal{N}(0,V_W(\phi,\Lambda)\int_{s_0}^1 s^{1-d}
\,ds)$, a normal distribution with mean~$0$ and variance $V_W(\phi,\Lambda
)\int_{s_0}^1 s^{1-d} \,ds$.
\end{theorem}

We turn now to the functional convergence of the total length process in
the planar case. Write
\[ 
\tau(s,R) := \exp\bigl([\log R - \log\log R](s-1)\bigr) = R^{s-1} (\log
R)^{1-s},
\]
and define the \textit{total} length process
\[ 
\nonumber\mathcal{L}^{R,W}_s := \frac{1}{R \sqrt{\log R}} \bar\Sigma_{\phi
}(Y(\tau(s,R),W_R)),\qquad s \in[0,1].
\]

\begin{theorem}\label{CLT2D}
The total length process $(\mathcal{L}^{R,W}_s)_{s\in[0,1]}$ converges in
law, as $R\to\infty,$ on the space
$\mathcal{D}[0,1]$ of c\`adl\`ag functions on $[0,1]$, endowed with the
Skorokhod $J_1$-topology, to\vadjust{\goodbreak}
$(\sqrt{V_W(\phi,\Lambda)} \mathcal{W}_s)_{s\in[0,1]}$, where, again,
$\mathcal{W}_{(\cdot)}$ stands for the
standard Wiener process. In particular, $(R\sqrt{\log R})^{-1}\bar
\Sigma_\phi(Y(1,W_R))$ converges in law to $\mathcal{N}(0,V_W(\phi,\Lambda
))$ with, again, $V_W(\phi,\Lambda)$ given by (\ref{VPHI}) or~(\ref
{VPHIALTERNATIVE}).
\end{theorem}

\begin{remark} We consider in this paper facet functionals of the
form $\phi(\cdot)=\Vol_{d-1}(\cdot)\zeta(\vec{\mathbf{n}}(\cdot))$; see
(\ref{PHIFORM}). Taking $\zeta\equiv1$, Theorems \ref
{BrownianConvergence}--\ref{CLT2D} reduce to the total surface area
case discussed in Theorems~\ref{limit1} and~\ref{limit2} in the
\hyperref[sec1]{Introduction}. However, the additional flexibility induced by the
introduction of $\zeta(\cdot)$ implies that our results allow us to
conclude limit theorems that are sensitive with respect to direction.
Taking, for example, $\zeta(\cdot)=\mathbf{1}\{\vec{\mathbf{n}}(\cdot)\in
U(n)\}$ to be the indicator function of a small neighborhood $U(n)$ of
a fixed direction $n\in\mathcal{S}_{d-1}$ satisfying $\mathcal{R}(U(n))>0$
[recall the decomposition (\ref{LADEF})], yields central limit theorems
also for the collection of tessellation facets having their normals in
$U(n)$. This means that our results are not only valid for the whole
STIT tessellation, but also for parts in arbitrary space directions.
\end{remark}

\begin{remark} So far we have restricted our considerations to space
dimensions $d\geq2$. STIT tessellations and their limit theory on the
line can also be considered. However, in~\cite{NW03} it was shown that
a STIT tessellation on $\mathbb{R}$ is nothing than a homogeneous Poisson
point process, or more precisely the intervals between its points.
These point processes and their limit theory are well known, and for
this reason we have focused on the cases $d\geq2$.
\end{remark}

\section{Proofs}\label{secPROOFS}

After having rephrased some background material on STIT tessellations
in the previous section, we are now prepared to present the proofs of
our limit theorems. Let us briefly recall that we will deal with a
fixed translation-invariant hyperplane measure $\Lambda$, and for this
reason we shall write, for example, $Y(t)$ instead of $Y(t\Lambda)$
without confusion. Moreover, we fix some compact and convex set
$W\subset\mathbb{R}^d$ having interior points and write $W_R=RW$ for $W$
dilated by a factor $R>0$. Moreover, recall that the facet functionals
we are dealing with have the representation (\ref{PHIFORM}), that
$\Sigma_\phi(Y(t))$ was defined in (\ref{FDEF}) and that the bar
notation $\bar\Sigma_\phi(Y(t))$ stands for the centered version $\Sigma
_\phi(Y(t))-\mathbb{E}\Sigma_\phi(Y(t))$.
\begin{pf*}{Proof of Theorem \protect\ref{BrownianConvergence}}
Notice first that because of $\Lambda([W_R])=R\Lambda([W])$, we have
%
\begin{eqnarray}\label{NormVar}
& & \frac{1}{R^d} A_{\phi^2}(Y(1,W_R))\nonumber\\
&&\qquad= \frac{1}{R} \int_{[W_R]} \frac{1}{R^{d-1}}
\zeta^2(\vec\mathbf{n}(H)) \sum_{f \in\Cells(Y(1,W_R)\cap H)} \Vol
^2_{d-1}(f) \Lambda(dH)\\
&&\qquad= \int_{[W]} \frac{1}{R^{d-1}}
\zeta^2(\vec\mathbf{n}(H)) \sum_{f \in\Cells( Y(1,W_R)\cap RH)} \Vol
^2_{d-1}(f) \Lambda(dH).\nonumber
\end{eqnarray}
We claim that, upon letting $R\to\infty$, this converges in probability to
\begin{eqnarray}\label{VPHI}
&& V_W(\phi,\Lambda)\nonumber\\
&&\qquad :=
\int_{[W]} \zeta^2(\vec\mathbf{n}(H)) \Vol_{d-1}(W \cap H)
\nonumber
\\[-8pt]
\\[-8pt]
\nonumber
&&\qquad{}\hspace*{32pt}\times  \frac{\mathbb{E}\Vol^2_{d-1}(\TypicalCell(Y(1) \cap H))}
{\mathbb{E}\Vol_{d-1}(\TypicalCell(Y(1)\cap H))} \Lambda(dH)
\\
&&\qquad= \Vol_d(W) \int_{\mathcal{S}_{d-1}} \zeta^2(u) \frac{\mathbb{E}\Vol
^2_{d-1}(\TypicalCell(Y(1) \cap u^{\bot}))}
{\mathbb{E}\Vol_{d-1}(\TypicalCell(Y(1) \cap u^{\bot}))}
\mathcal{R}(du),\nonumber
\end{eqnarray}
where $\mathcal{R}$ is the directional distribution of the tessellation as
given by (\ref{LADEF}). To see it, recall that $ Y(1)\cap RH$ is a STIT
tessellation in $RH$ for each $R >0$ and $H \in\mathcal{H}.$ Thus,
applying \cite[(4.6) and Theorem 4.1.3]{SW} to this tessellation and
the fact that $\mathbb{E}\Vol_d(\TypicalCell(Y(1)\cap u^\perp))$ is the
same as the inverse cell density of the tessellation $Y(1)\cap u^\perp$
[see (10.4) ibidem], we get
%
\begin{eqnarray}\label{EXPE}
& & \lim_{R\to\infty} \frac{1}{R^{d-1}} \mathbb{E} \sum_{f \in
\Cells(Y(1,W_R) \cap RH)} \Vol^2_{d-1}(f)
\nonumber
\\[-8pt]
\\[-8pt]
\nonumber
&&\qquad= \Vol_{d-1}(W \cap H) \frac{\mathbb{E}\Vol^2_{d-1}(\TypicalCell(Y(1)
\cap H))}
{\mathbb{E}\Vol_{d-1}(\TypicalCell(Y(1) \cap H))}.
\end{eqnarray}
Next, we observe that $Y(1,W_R) \cap R H \stackrel{D}{=} Y(1)\cap R \cdot
_H (H \cap W)$, where
$\cdot_H$ is the scalar multiplication relative in $H$, that is, to say,
$H \ni R \cdot_H x = p_H(0) + R (x-p_H(0)),  x \in H$
with $p_H$ standing for the orthogonal projection on $H.$ Thus, using
the recently developed strong mixing and tail triviality theory for
STIT tessellations, especially \cite[Theorem 2]{LR}, noting that tail
trivial stationary processes are ergodic \cite[Proposition 14.9]{GEO}
and then applying the multidimensional ergodic theorem (see, e.g.,
Corollary 14.A5 ibidem),
to $\frac{1}{R^{d-1}} \sum_{f \in\Cells(Y(1)\cap R \cdot_H (H \cap W)
)} \Vol^2_{d-1}(f),$ we get from (\ref{EXPE}) that
\begin{eqnarray*}
& & \lim_{R\to\infty} \frac{1}{R^{d-1}} \sum_{f \in\Cells
(Y(1,W_R)\cap R H)} \Vol^2_{d-1}(f)\\
 &&\qquad= \Vol_{d-1}(W \cap H) \frac{\mathbb{E}\Vol
^2_{d-1}(\TypicalCell(Y(1)\cap H))}
{\mathbb{E}\Vol_{d-1}(\TypicalCell(Y(1)\cap H))}
\end{eqnarray*}
in probability. Putting this together with (\ref{NormVar}) and
integrating over $[W]$ yields
%
\begin{equation}\label{PRZEDZBIEZN}
\lim_{R \to\infty} \frac{1}{R^d} A_{\phi^2}(Y(1,W_R)) = V_W(\phi
,\Lambda)\qquad \mbox{in probability},
\end{equation}
as required.

Note now that by the scaling properties of $Y(s,W_R)$ and $\phi^2$ for
$s > 0$, we have
%
\begin{eqnarray}\label{SCALINGIDS}
\frac{1}{R^d} A_{\phi^2}(Y(s,W_R)) &\stackrel{D}{=} &\frac
{1}{R^d s^{2d-1}} A_{\phi^2}(Y(1,W_{sR}))
\nonumber
\\[-8pt]
\\[-8pt]
\nonumber
&\stackrel{D}{=}&
\frac{1}{s^{d-1}} \frac{1}{(Rs)^d} A_{\phi^2}(Y(1,W_{sR})).
\end{eqnarray}
Thus, combining (\ref{PRZEDZBIEZN}) with the scaling relation (\ref
{SCALINGIDS}), we get
%
\begin{equation}\label{CONVERGSTAT}
\lim_{R\to\infty} \frac{1}{R^d} A_{\phi^2}(Y(s,W_R)) = \frac
{1}{s^{d-1}} V_W(\phi,\Lambda)
\end{equation}
in probability and uniformly in $s \in[s_0,1].$

This crucial statement puts us now in context of general martingale
limit theory. Indeed, using Propositon~\ref{propMART3}, we
see that
\[
\mathcal{S}^{R,W}_{s_0,s} = \frac{1}{R^{d/2}} [\bar\Sigma_{\phi
}(Y(1,W_R)) - \bar\Sigma_{\phi}(Y(s_0,W_R))]
\]
is a martingale with absolutely continuous predictable quadratic
variation process
%
\begin{equation}\label{QUADVAR}
\langle\mathcal{S}^{R,W}_{s_0,\cdot} \rangle_t = \int_{s_0}^t \frac
{1}{R^d} A_{\phi^2}(Y(s,W_R)) \,ds
\end{equation}
by (\ref{MART3}) and (\ref{PRZEDZBIEZN}); see again Paragraph I.8 in
\cite{LS}.
In these terms, (\ref{CONVERGSTAT}) yields for each~$t$,
%
\begin{equation}\label{QVARCONVREL}
\lim_{R \to\infty} \langle\mathcal{S}^{R,W}_{s_0,\cdot} \rangle_t = \int
_{s_0}^t \frac{1}{s^{d-1}} V_W(\phi,\Lambda)\,ds\qquad \mbox{in
probability. }
\end{equation}
We now want to apply the martingale functional central limit theorem.
Whereas this is well known for continuous martingales, we need a
version for martingales in the Skorokhod space $\mathcal{D}[s_0,1]$. In
this paper, we will make use of the version formulated as Theorem 2.1
in the survey article~\cite{WW}. In order to apply this theorem,
several conditions have to be checked. Condition
(ii.6) in [\cite*{WW}, Theorem 2.1] is just~(\ref{QVARCONVREL}), whereas
condition (ii.4) there is trivially verified, because the predictable
quadratic variation $\langle\mathcal{S}^{R,W}_{s_0,\cdot} \rangle$ has no
jumps by (\ref{QUADVAR}).
It remains to check condition~(ii.5) ibidem, which is that the second
moment of the maximum jump
$\mathcal{J}(\mathcal{S}^{R,W}_{s_0,\cdot};1)$ of the process $(\mathcal{S}^{R,W}_{s_0,s})_{s \in[s_0,1]}$
goes to $0$, as $R \to\infty.$ More precisely,
\[
\mathcal{J}(\mathcal{S}^{R,W}_{s_0,\cdot};1)=\sup_{s_0\leq t\leq1}|\mathcal{S}^{R,W}_{s_0,t}-\mathcal{S}^{R,W}_{s_0,t-}|,
\qquad  \mathcal{S}^{R,W}_{s_0,t-}=\lim_{s\uparrow t}\mathcal{S}^{R,W}_{s_0,s},
\]
and we have to check that
\[
\lim_{R\rightarrow\infty}\mathbb{E}\mathcal{J}^2(\mathcal{S}^{R,W}_{s_0,\cdot};1)=0.
\]
To this end, note first that, with probability one,
$\mathcal{J}(\mathcal{S}^{R,W}_{s_0,\cdot};1)$ is bounded from above by a
constant multiple of
$R^{-d/2}$ times the $(d-1)$th
power of the diameter of the largest cell of $Y(s_0,W_R).$ Since the
typical cell distribution of $Y(s_0)$
is the same as that of a Poisson hyperplane tessellation with
intensity measure $s_0 \Lambda$ (see Theorem 1 in~\cite{STP1} or
Section~\ref{secSTIT} above), we conclude by standard properties of
Poisson hyperplane tessellations that the expected number
$e(Y(s_0,W_R),D)$ of cells in $Y(s_0,W_R)$ with diameter exceeding $D$
is of order $O(R^d \exp(-D))$. To see it, write $\diam(c)$ for the
diameter of a cell $c$ and $\mathbf{1}\{\cdot\}$ for the usual indicator
function and rewrite $e(Y(s_0,W_R),D)$ as
\[
e(Y(s_0,W_R),D)=\mathbb{E}\sum_{c\in\Cells(Y(s_0,W_R))}\mathbf{1}\{\diam
(c)>D\},
\]
which by Theorem 4.1.3 in~\cite{SW} is of the same order as the mean
number $N(Y(s_0,W_R))=\mathbb{E}\sum_{c\in\Cells(Y(s_0,W_R))}1$ of cells
in $Y(s_0,W_R)$ times the probability that the typical cell diameter of
$Y(s_0)$ exceeds $D$ (the additional condition in Theorem 4.3.1 in~\cite
{SW} is easily verified by using Steiner's formula together with the
fact that the typical cell of a Poisson hyperplane tessellation has
finite mean intrinsic volumes; see Theorem 10.3.3 ibidem). Thus,
$e(Y(s_0,W_R),D)$ satisfies
\[
e(Y(s_0,W_R),D)=O\bigl(N(Y(s_0,W_R))\mathbb{P}\bigl(\diam(\TypicalCell
(Y(s_0)))>D\bigr)\bigr).
\]
Stationarity of the tessellation implies that the first factor is of
volume order $R^d$. We claim that the second term is bounded from above
by $c_1e^{-c_2s_0 D}$, where $c_1$ and $c_2$ are constants which depend
on the hyperplane measure $\Lambda$. To this end we notice first that
the typical cell of $Y(s_0)$ is stochastically smaller than the almost
surely uniquely determined cell $Z_0$ of $Y(s_0)$ containing the
origin; cf. Corollary 10.4.1 in~\cite{SW}. Moreover, equation (20) in
\cite{HS1} with $\Sigma=\diam$, $\varepsilon=0$ and $\alpha=1$ there
(the other parameters are then $k=\tau=1$) implies that there are
constants $c_1,c_2>0$ depending on the hyperplane measure $\Lambda$
such that $\mathbb{P}(\diam(Z_0)>D)\leq c_1e^{-c_2s_0D}$. This implies
\[
\mathbb{P}\bigl(\diam\bigl(\TypicalCell(Y(s_0))\bigr)>D\bigr)\leq\mathbb{P}\bigl(\diam(Z_0)>D\bigr)\leq
c_1e^{-c_2s_0D}.
\]
Putting these two issues together leads to the desired order $O(R^d\exp
(-D))$ for the expected number of cells in $Y(s_0,W_R)$ with diameter
exceeding $D$. Recalling that $\mathcal{J}(\mathcal{S}^{R,W}_{s_0,\cdot};1)$
is bounded from above by a constant multiple of
$R^{-d/2}$ times the $(d-1)$th power of the diameter of the largest
cell of $Y(s_0,W_R)$, and putting $u = D^{d-1} R^{-d/2}$ we find
%
\begin{equation}\label{JBOUND}
\mathbb{P}\bigl(\mathcal{J}(\mathcal{S}^{R,W}_{s_0,\cdot};1) > u\bigr) = O\bigl(R^d \exp
\bigl(-R^{d/(2d-2)} u^{1/(d-1)}\bigr)\bigr).
\end{equation}
Clearly, (\ref{JBOUND}) is sufficient to guarantee that
\[
\lim_{R\to\infty} \mathbb{E}\mathcal{J}^2(\mathcal{S}^{R,W}_{s_0,\cdot};1) = 0,
\]
which gives the required condition (ii.5) of Theorem 2.1 in~\cite{WW}.
This result yields now the functional convergence as stated in our
Theorem~\ref{BrownianConvergence}.
\end{pf*}

Before turning to the proof of Theorem~\ref{CLT2D} we provide an
alternative formula for the factor $V_W(\phi,\Lambda)$. Readers not
specialized in convex or stochastic geometry could also skip this
alternative representation and directly jump to the next paragraph,
because Proposition~\ref{propVW} will not be used in the sequel.
However, having such a more explicit variance expression is useful for
other purposes and has already been used in our work~\cite{ST3}. We
denote, as in~\cite{SW} or~\cite{STP1}, by $\Pi$ the associated zonoid
of a Poisson hyperplane tessellation with intensity measure $\Lambda$,
by $\Pi^o$ its polar body and by $\mathcal{R}$ the directional
distribution of the STIT tessellation from (\ref{LADEF}); see~\cite{SW}
for the precise definitions of $\Pi$ and $\Pi^o$.
\begin{proposition}\label{propVW} We have
%
\begin{equation}V_W(\phi,\Lambda)=\Vol_{d}(W){(d-1)!\over2^{d-1}}\int
_{\mathcal{S}_{d-1}}\zeta^2(u)\Vol_{d-1}((\Pi|u^\perp)^o)\mathcal{R}(du),\label{VPHIALTERNATIVE}
\end{equation}
where $\Pi|u^\perp$ stands for the orthogonal projection of $\Pi$ onto
the hyperplane $u^\perp$, and where the
polar body $(\Pi|u^{\perp})^o$ is considered relative to $u^{\perp}.$
In the isotropic case, that is, when $\mathcal{R}=\nu_{d-1}$, the uniform
distribution on the unit sphere $\mathcal{S}_{d-1}$, this reduces to
\[
V_W(\phi,\Lambda_\mathrm{iso})=\Vol_d(W)2^{d-1}\pi^{d-{3/2}}{\Gamma
({(d+1)/2})^{d-1}\over\Gamma({d/2}
)^{d-2}}\int_{\mathcal{S}_{d-1}}\zeta^2(u)\nu_{d-1}(du).
\]
In particular for $\zeta\equiv1$, $W=B_1^d$ the unit ball and $d=2$
and $d=3$, we conclude the exact values
\[
V_{B_1^2}(\Vol_1,\Lambda_\mathrm{iso})=\pi^2  \quad   \mbox{and} \quad   V_{B_1^3}(\Vol_2,\Lambda_\mathrm{iso})={32\over3}\pi^2.
\]
\end{proposition}
\begin{pf}
At first,~\cite{FW}, Corollary 3.7, provides a general formula for the second moment of
the volume of the typical Poisson cell of a stationary Poisson
hyperplane tessellation $\PHT(\Lambda)$ in $\mathbb{R}^d$ having
intensity measure~$\Lambda$. In terms of the zonoid $\Pi$ it reads
\[
\mathbb{E}\Vol_d^2(\TypicalCell(\PHT(\Lambda)))={d!\over2^d}{\Vol_d(\Pi
^o)\over\Vol_d(\Pi)},
\]
where we have used [\cite*{SW}, equation (4.63)]. Moreover, the first
volume moment of the typical cell of a Poisson hyperplane tessellation
with intensity measure $\Lambda$ equals $1/\Vol_d(\Pi)$ according to
[\cite*{SW}, Theorem 10.3.3 and (10.4)]. Using now equation (4.61) ibidem
and the fact that STIT tessellations have Poisson typical cell
distributions and replacing $d$ by $d-1$ in the last two formulas, we
obtain (\ref{VPHIALTERNATIVE}) immediately from (\ref{VPHI}). The
precise value in the stationary and isotropic case can be calculated
from the fact that in this case, $\Pi$ is a ball with a known radius;
see~\cite{SW}.\vadjust{\goodbreak}
\end{pf}

\begin{pf*}{Proof of Theorem \protect\ref{CLT2D}}
Recall that $\tau(s,R)$ is defined by
\[
\tau(s,R)=\exp\bigl([\log R-\log\log R](s-1)\bigr)=R^{s-1}(\log R)^{1-s},
\]
and note that this implies
%
\begin{eqnarray}\label{TAUPROP}
\tau(0,R)& =& \frac{\log R}{R},\qquad  \tau(1,R) = 1,
\nonumber
\\[-8pt]
\\[-8pt]
\nonumber
 \frac{\partial
}{\partial s}\tau(s,R)& = &\tau(s,R) [\log R - \log\log R].
\end{eqnarray}
Thus, defining the auxiliary process
\begin{eqnarray*}
 M^{R,W}_s &= &M_s \\
&:=&
\frac{1}{R \sqrt{\log R - \log\log R}} [\bar\Sigma_{\phi}(Y(\tau
(s,R),W_R)) - \bar\Sigma_{\phi}(Y(\tau(0,R),W_R))]
\end{eqnarray*}
and using (\ref{MART3}) with $W_R := R W$ and under variable
substitution $s := \tau(u,R)$ and $t := s,$
with left-hand side variables corresponding to the notation of (\ref
{MART3}), and those on the right-hand side to that used here,
we see that, by (\ref{TAUPROP}),
\[ 
(M_s)_{s\in[0,1]}\quad \mbox{and} \quad\biggl( M^2_s - \int_0^s
\frac{\tau(u,R)}{R^2} A_{\phi^2}(Y(\tau(u,R), W_R)) \,du \biggr)_{s\in[0,1]}
\]
are $\Im_{\tau(s,R)}$-martingales. In particular (see, once more, \cite
[Paragraph I.8]{LS}), the predictable quadratic variation process
$\langle M \rangle_s$ is given by
%
\begin{equation}\label{MQV}
\langle M \rangle_s = \int_0^s \frac{\tau(u,R)}{R^2} A_{\phi^2}(Y(\tau
(u,R),W_R)) \,du, \qquad s \in[0,1].
\end{equation}
Repeating the argument leading to (\ref{CONVERGSTAT}) we see that
%
\begin{equation}\label{CONVSTAT2}
\lim_{R\to\infty} \frac{\tau(s,R)}{R^2} A_{\phi^2}(Y(\tau(s,R), W_R))
= V_W(\phi,\Lambda)
\end{equation}
in probability and uniformly in $s \in[0,1].$ Note that the uniformity in
$s$ comes, as in the case of (\ref{CONVERGSTAT}), from the relation
(\ref{SCALINGIDS})
implying that, in distribution, all instances of the left-hand side for
different values of $s$ are just re-scalings of the same
object $\tilde{R}^{-2} A_{\phi^2}(Y(1,W_{\tilde{R}}))$ for $\tilde{R} =
R / \tau(s,R)$, and thus, in terms of the considered
convergence in probability to a deterministic limit, we are just
dealing with a single asymptotic statement. Consequently, by (\ref
{CONVSTAT2}) and in full analogy to (\ref{QVARCONVREL}),
%
\begin{equation}\label{QVARCONVREL2}
\lim_{R\to\infty} \langle M \rangle_s = \int_0^s V_W(\phi,\Lambda) \,du
= s V_W(\phi,\Lambda)\qquad \mbox{in probability.}
\end{equation}
Thus, we are again in a position to apply the martingale functional
central limit theorem \cite[Theorem 2.1]{WW} yielding the functional
convergence in law, as $R \to\infty,$ in $\mathcal{D}[0,1]$ of
$(M_s)_{s\in[0,1]}$ to the random process $(\sqrt{V_W(\phi,\Lambda)}
\mathcal{W}_s)_{s\in[0,1]}.$ Indeed,
condition (ii.6) there is just (\ref{QVARCONVREL2}), condition (ii.4)
is trivial in view of (\ref{MQV}), whereas
the condition (ii.5) is verified by noting that, with\vadjust{\goodbreak} probability one,
$\mathcal{J}(M;1) = \frac{1}{R \sqrt{\log R}}
O(R \diam(W)) = O(1/\sqrt{\log R})$, so that in particular
\[
 \lim_{R\to\infty} \mathbb{E}\mathcal{J}^2(M;1) = 0,
\]
as required.

Denoting now by
$C^{R,W}$ the \textit{correction term} ${1\over R\sqrt{\log R}}\bar\Sigma
_{\phi}(Y(\tau(0,R),W_R)),$ such that
\[
\mathcal{L}^{R,W}_s = C^{R,W} + \sqrt{\frac{\log R - \log\log R}{\log R}}
M_s,
\]
noting that $\log R - \log\log R \sim\log R$ and that, by the scaling
property of STIT tessellations and
by (\ref{EQVARASYM2D}),
\begin{eqnarray*}\label{NEGcorr}
 \Var(C^{W,R}) &=& O\bigl([R^{-2} (\log R)^{-1}] [R^2/(\log
R)^2] [(\log R)^2 (\log\log R)]\bigr)\\
&=& O(\log\log R/ \log R),
\end{eqnarray*}
we see that the processes $M_s$ and $\mathcal{L}^{R,W}_s$ are
asymptotically equivalent in $\mathcal{D}[0,1]$, as $R\to\infty.$
This completes the proof of Theorem~\ref{CLT2D}.
\end{pf*}

\begin{remark}
In the context of proof of Theorem~\ref{CLT2D} it should be remarked
that the ``negligible correction term''
$C^{R,W}$ has its variance of order
\[
O(\log\log R / \log R)
\]
and is thus indeed tending to $0,$
but extremely slowly. Consequently, although the Gaussian CLT holds
for $\mathcal{L}^{R,W}_1,$ it is quite natural
to expect that the convergence rates are extremely slow and
conjecturedly logarithmic. This is due to
the fact that dimension $2$ is the largest dimension (critical
dimension) where the Gaussian limits are still
present. In dimensions $3$ and higher there is no Gaussian CLT and the
``correction term'' analogous to
$C^{R,W}$ will turn out order-determining rather than negligible, as
shown by Theorem~\ref{thmNONGAUSSIAN}.
\end{remark}

We turn now to the higher-dimensional situations. Even if in the
formulation of Theorem~\ref{thmNONGAUSSIAN} we have used the surface
area functional, we will show the statement in a more general context,
where ${\Vol}_{d-1}(Y(1,W_R))$ is replaced by a general cumulative
facet functional $\Sigma_\phi(Y(1,W_R))$ satisfying (\ref{PHIFORM}).

We claim that the argument from the proof of Theorem~\ref{CLT2D} cannot
be repeated for $d > 2.$ Intuitively, this is due to the fact that for
$d > 2$ the variance order of $\bar\Sigma_{\phi}(Y(1,W_R))$ is
$O(R^{2(d-1)})$, see below, whereas the variance order of the increment
$\bar\Sigma_{\phi}(Y(1,W_R))-\bar\Sigma_{\phi}(Y(s_0,W_R))$, with some
time instant $0<s_0<1$, is $O(R^d)$ as seen from Theorem \ref
{BrownianConvergence}. Hence, for $d>2$ we conclude that even the very
first facets born in the MNW-cell-division process already bring a
nonnegligible contribution to the overall variance. Thus, we cannot
split the whole STIT construction into the \textit{warm-up phase} ($t \in
[0,R^{-1} \log R]$ for $d=2$) with negligible variance\vadjust{\goodbreak} contribution and
the \textit{proper phase} unfolding already in a typical STIT environment.
In fact, the claim is that the CLT does not hold for STIT surface
functionals in dimension greater than $2!$

\begin{pf*}{Proof of Theorem~\ref{thmNONGAUSSIAN}}
Recall that we do not show this fact in full generality for all
nondegenerate hyperplane measures $\Lambda$ and all windows $W$, but
restrict ourself to a particular case, where $\Lambda$ is given by (\ref
{PARTICULARLA}) and where $W=[0,1]^d$. To see the non-Gaussianity,
observe first that, by the scaling property of STIT tessellations,
%
\begin{equation}\label{RNIEDOCTG}
R^{-(d-1)} \bar\Sigma_{\phi}(Y(1,W_R)) \stackrel{D}{=} \bar\Sigma_{\phi}(Y(R,W)),
\end{equation}
which implies that the variance $\Var(\Sigma_{\phi}(Y(1,W_R)))$ is of
order $O(R^{2(d-1)})$. Indeed, this follows directly from the special
form (\ref{PHIFORM}) of the facet functional $\phi$ and the scaling
relation $Y(1,W_R)\stackrel{D}{=}RY(R,W)$. Further, recall that by (\ref
{MART3}) the process $R \mapsto\bar\Sigma_{\phi}(Y(R,W))$ is a
square-integrable
martingale with absolutely continuous predictable quadratic variation
process given as in~(\ref{PREDQUVAR}).
Moreover, by Proposition~\ref{thmvariance} we conclude
\begin{eqnarray*}
&&\Var\bigl(\Sigma_{\phi}(Y(R,W))\bigr)\\
&&\qquad =
\int_{[W]} \zeta^2(\vec\mathbf{n}(H)) \int_{W \cap H}
\int_{W \cap H} \frac{1-\exp(-R\Lambda([xy]))}{\Lambda([xy])} \,dx \,dy
\Lambda(dH),
\end{eqnarray*}
which is bounded uniformly in $R.$
Consequently, by the martingale convergence theorem (cf. Corollary 7.22
in~\cite{K}) there exists a centered square-integrable
random variable $\Xi(W)$ such that
\[ 
\Xi(W) = \lim_{R\to\infty} \bar\Sigma_{\phi}(Y(R,W))
\]
a.s. and in $L^2$ and, moreover,
%
\begin{eqnarray}\label{XIVAR}
\Var(\Xi(W)) &=& \lim_{R\rightarrow\infty}\Var\bigl(\Sigma_\phi
(Y(R,W))\bigr)
\nonumber
\\[-8pt]
\\[-8pt]
\nonumber
&=& \int_{[W]} \zeta^2(\vec\mathbf{n}(H)) \int_{W \cap H}
\int_{W \cap H} \frac{1}{\Lambda([xy])} \,dx \,dy \Lambda(dH).
\end{eqnarray}
Using now (\ref{RNIEDOCTG}) we readily conclude that
\[ 
 R^{-(d-1)} \bar\Sigma_{\phi}(Y(1,W)) \Longrightarrow\Xi(W),
\]
as $R \to\infty,$ where $\Longrightarrow$ stands for convergence in
distribution.

We show now that the variable $\Xi(W)$ cannot be Gaussian. To see it,
consider the event $\mathcal{E}_{N},  N >0,$ that only hyperplanes orthogonal
to $e_1$ have been born during time $[0,1]$ in the MNW-construction and
that their number exceeds
$N.$ Observe that, in view of the special form (\ref{PARTICULARLA}) of
the hyperplane measure $\Lambda,$
$ \mathbb{P}(\mathcal{E}_N) = \exp(-d) \sum_{k=N+1}^{\infty} \frac{1}{k!} $
and thus
%
\begin{equation}\label{TAILSN}
\mathbb{P}(\mathcal{E}_N) = \exp(- \Theta(N \log N)),
\end{equation}
where by $\Theta(\cdot)$ we mean a function bounded both from below and
above by multiples of the argument. Further, given the fixed collection\vadjust{\goodbreak}
of all hyperplanes $H_1,\ldots,H_k$ ($k>N$) born at times between $0$
and $1,$ on the event $\mathcal{E}_N$, we see that the conditional law of
$\Xi(W)$
coincides with that of $k-d$ plus the sum of independent copies $\xi
_1,\ldots,\xi_{k+1}$ of
$\Xi(W_1),\ldots,\Xi(W_{k+1})$ respectively,
where $W_j,  j=1,\ldots,k+1,$ are the parallelepipeds into which
$W=[0,1]^d$ is partitioned by $H_1,\ldots,H_k$. More formally, we have
the relation
\[
\mathbb{P}\bigl(\Xi(W)>u|\mathcal{E}_{N,k}\bigr)=\mathbb{P}\bigl(k-d+(\xi_1+\cdots+\xi
_{k+1})>u\bigr),   \qquad   u\in\mathbb{R},
\]
where $\mathcal{E}_{N,k}$ is the event that exactly the hyperplanes
$H_1,\ldots,H_k$ ($k>N$) orthogonal to $e_1$ are born within time
$[0,1]$. Note that the extra $k$ above is the sum of the
$(d-1)$-volumes of $W\cap H_i$, whereas $-d = - \mathbb{E}\Sigma_{\phi
}(Y(1,W))$ is the centering term.

Since $\Var(\xi_1+\cdots+\xi_{k+1}) =
\sum_{j=1}^{k+1} \Var(\Xi(W_j))$, which is bounded from above by $\Var
(\Xi(W))$ in view of (\ref{XIVAR}),
by Chebyshev's inequality we get
\begin{eqnarray*}
\mathbb{P}\bigl(\xi_1+\cdots+\xi_{k+1} \geq- 2 \sqrt{\Var(\Xi(W))}\bigr)
&\geq& 1-{\Var(\xi_1+\cdots+\xi_{k+1})\over4\Var(\Xi(W))}\\
&\geq& 1-{1\over4}={3\over4}.
\end{eqnarray*}
Thus, in view of (\ref{TAILSN}) we end up with
\[ 
\mathbb{P}\bigl(\Xi(W) > N\bigr) \geq\frac{3}{4} \mathbb{P}\bigl(\mathcal{E}_{N + 2\sqrt{\Var(\Xi(W))} + d}\bigr)
= \exp(-\Theta(N\log N)).
\]
Since Gaussian variables exhibit tail decay of the order $\exp(-\Theta
(N^2)),$ the random variable
$\Xi(W)$ cannot be Gaussian.
\end{pf*}

\section*{Acknowledgments}
The second author is indebted to Werner Nagel (Jena), Matthias Reitzner
(Osnabr\"uck) and Pierre Calka (Rouen) for interesting discussions,
suggestions, hints and remarks. Moreover, the comments of an anonymous
referee were helpful in improving the presentation and the style of the paper.



\printaddresses

\end{document}